\documentclass[final]{siamltex}
\usepackage{amsfonts}
\usepackage{amssymb,amsmath}
\usepackage[mathscr]{euscript}

\DeclareMathAlphabet{\itbf}{OML}{cmm}{b}{it}
 \DeclareMathAlphabet\mathbfcal{OMS}{cmsy}{b}{n}
 
\def\EE{\mathbb{E}}

\def\RR{\mathbb{R}}
\def\NN{\mathbb{N}}
\def\eps{\varepsilon}

\def\bq{{{\itbf q}}}

\def\by{{{\itbf y}}}
\def\bx{{{\itbf x}}}

\def\br{{{\itbf r}}}
\def\bX{{{\itbf X}}}

\newcommand{\bka}{{\boldsymbol{\kappa}}}

\newcommand{\la}{\lambda}

\newcommand{\vx}{\vec{\itbf x}}
\newcommand{\vy}{\vec{\itbf y}}
\newcommand{\vecr}{\vec{\itbf r}}
\newcommand{\vna}{\vec{\nabla}}
\newcommand{\om}{\omega}

\newcommand{\cR}{{\mathscr R}}
\newcommand{\ep}{\varepsilon}
\newcommand{\what}{\widehat}
\newcommand{\bphi}{\boldsymbol{\phi}}
\newcommand{\cS}{{\mathcal{S}}}

\newcommand{\vka}{\vec{\bka}}
\newcommand{\vk}{\vec{\boldsymbol{\mathscr{K}}}}

\newcommand{\mk}{\mathscr{K}}
\newcommand{\bmk}{\boldsymbol{\mk}}

\begin{document}

\title{Derivation of a one-way radiative transfer equation in  random media}
\author{Liliana Borcea \footnotemark[1] and Josselin Garnier \footnotemark[2]}

\footnotetext[1]{Department of Mathematics,
University of Michigan,
Ann Arbor, MI 48109-1043
borcea@umich.edu
}
\footnotetext[2]{Laboratoire de Probabilit\'es et Mod\`eles Al\'eatoires
\& Laboratoire Jacques-Louis Lions,
Universit{\'e} Paris Diderot,
75205 Paris Cedex 13,
France. garnier@math.univ-paris-diderot.fr}
\maketitle

\begin{abstract}
We derive from first principles a one-way radiative transfer equation
for the wave intensity resolved over directions (Wigner transform of the wave field) in random media. It is an initial value 
problem  with excitation from a source which emits waves in a preferred,
forward direction. The equation is derived in a regime with small random fluctuations of the wave speed but 
long distances of propagation with respect to the wavelength, so that cumulative scattering is
significant. The correlation length of the medium and the scale of the support of the source 
are slightly larger than the wavelength, and the waves propagate in a wide cone with opening angle less than $180^o$, so that  the backward and evanescent waves are negligible. 
The scattering regime is a bridge between that of radiative transfer, where the waves propagate in 
all directions and the paraxial regime, where the waves propagate in a narrow angular cone. We connect
the one-way radiative transport equation with the equations satisfied by the Wigner transform of the wave field in these regimes.
\end{abstract}

\section{Introduction}
Light propagation in scattering media can be modeled by a boundary value problem for the radiative transfer equation \cite{ishimaru,ryzhik96, bal2008transport}.  
The light intensity resolved over directions, also known as the Wigner transform of the wave field, satisfies this equation 
with incoming boundary conditions on the illuminated part of the boundary, and outgoing conditions on the remainder of the boundary.  
The problem is of interest  in applications such as optical tomography, where structural variations 
in tissue are to be determined from measurements of scattered light \cite{AS_REV}.

The derivation of the radiative transfer equation from the wave equation is a fundamental challenge.
Existing heuristic derivations  from the wave equation in random media, obtained 
when  the wavelength, the correlation length of the medium and the scale of variation of the source 
are of the same order, and much smaller than the propagation distance, use either 
multiscale asymptotic analysis  \cite{ryzhik96} or diagrammatic 
perturbation theory \cite{barabanenkov1968radiation,vollhardt1980diagrammatic}.
However, 
as discussed by Mandel and Wolf in their monography \cite{mandel}, or more recently in the tutorial \cite{schotland15},
there is no satisfactory  or rigorous derivation of the macroscopic theory of radiative transfer
from the microscopic theory of wave propagation in random media, except in some special cases.
Therefore, the rigorous derivation of a radiative transfer-like equation from the wave equation, beyond the special cases 
mentioned in these references, would be of interest for the radiative transfer community.

The radiative transfer equation  poses formidable computational challenges in optical tomography, where repeated solutions of the equation are needed to solve the inverse problem with optimization \cite{bal2008transport,AS_REV}. 
This is why a simplified  diffusion model is often used \cite{AS_REV},
where the medium is assumed optically thick, so that light is diffusive due to very strong scattering. 
This leads to considerable simplification, but may produce anomalies in the reconstructed images \cite{dehghani99}.
A recent study \cite{kim15} shows that in mesoscopic scattering regimes, where light penetrates to about 
one centimeter depth in tissue \cite{gonzalez2013robust}, scattering is forward-peaked and
a simpler one-way radiative transport model can be used, where the intensity satisfies an initial value problem.
The one-way radiative transfer equation is obtained in \cite{kim15}
from the standard radiative transfer equation
by simply ignoring the intensity in the backward 
directions.

Our first goal in this paper is to derive rigorously the one-way radiative transfer equation, from first principles, 
starting from the wave equation in random media. 
The second goal is to bridge between the mesoscopic scattering regime,  the standard radiative 
transfer regime on one side, and the paraxial approximation regime on the other side. We also connect to 
the diffusion approximation.

To derive the one-way transfer equation we consider waves in media with small 
random fluctuations of the wave speed, at long propagation distances with respect to the wavelength, 
where cumulative scattering effects are significant.  The typical size of the inhomogeneities, measured 
by the correlation length, and the scale of variation of the source are slightly larger than 
the wavelength, so that the waves propagate in an angular cone with axis along a preferred forward direction called range.
We analyze the propagation in this regime using a plane wave decomposition of the field, with amplitudes that 
are range dependent random fields. They satisfy a system of coupled stochastic differential equations driven by the 
random fluctuations of the wave speed, and can be analyzed in detail with probabilistic limit theorems. 
Consequently, we can quantify the loss of coherence of the wave field i.e., its randomization due to 
scattering, and derive the radiative transfer equation satisfied by the Wigner transform of the wave amplitudes.
The result extends the model proposed in \cite{kim15}, and defines the differential scattering cross-section 
and the total scattering cross-section in terms of the autocovariance of the fluctuations of the wave speed. 

Once we derive the one-way transfer equation we show that it is equivalent  to the standard radiative transfer equation \cite{ryzhik96} in regimes with negligible backscattering. We also connect to the diffusion approximation theory, by considering
the high-frequency limit of the equation. Transport in the paraxial approximation, which applies to waves propagating in a narrow angle cone
along the range axis,  is analyzed in \cite{garnier2009paraxial}, 
using the It\^o-Schr\"odinger model of wave propagation. Here we rediscover the results  starting from the 
one-way radiative transfer equation, in the high-frequency limit and for a large support of the source.

The paper is organized as follows: We begin in section 
\ref{sect:main} with the model of the random medium and the formulation of the problem.
The main results are stated in section \ref{sect:MR}. We describe the mean wave field and 
the randomization of the components of the wave quantified by the  scattering mean free paths.  
We also state the one-way radiative transfer equation. The connection to the 
equation in \cite{kim15} is in subsection \ref{sect:ResA}, and to the standard radiative transfer
theory in subsection \ref{ap:RT}. The connection to the paraxial approximation is in 
subsection \ref{app:C}.
The derivation of the results is  in section \ref{sect:deriv}. We begin with the scaling regime in subsection 
\ref{sect:multi}, and then give the wave decomposition in subsection \ref{sect:WD}.
The probabilistic limit of the wave amplitudes is studied in subsection \ref{sect:DL} and appendix 
\ref{ap:DL}. 
We use it to describe the evolution of the mean field in subsection \ref{sect:COH} and to derive 
the one-way radiative transfer equation for the Wigner transform
in subsection  \ref{sect:TEQ}. The high-frequency limit which leads to either the diffusion 
approximation or the paraxial approximation is studied in section \ref{sect:HF}. We end with a summary in 
section \ref{sect:sum}.

\section{Formulation of the problem}
\label{sect:main}%
The time-harmonic  field $u(\vx)$   satisfies the wave equation:
\begin{align}
\frac{\omega^2}{c^2(\vx)} u(\vx) +\Delta_{\vx} u(\vx) = -F\Big( \frac{\bx}{X} \Big) \delta(z) ,
  \label{eq:F0}
\end{align}
for $\vx =(\bx,z) \in \mathbb{R}^{d+1}$ and frequency $\om \in \mathbb{R}$. 
Here $\Delta_{\vx}$ is the Laplacian 
operator in $\mathbb{R}^{d+1}$ and since the frequency  is 
constant,  we suppress $\om$  from the arguments of $u$ and $F$. 
The excitation is due to a localized source $F$ which emits waves in the direction
$z$, called range.  The function $F$ depends on  the dimensionless vector $\br \in
\mathbb{R}^d$, and its magnitude  is negligible for $|\br| >
O(1)$, so that $X$ scales the spatial support of the source.  

The waves propagate in a linear medium with
speed of propagation $c(\vx)$ satisfying
\begin{equation}
\frac{1}{c^2(\vx)} = \frac{1}{c_o^2} \Big[1 + 1_{(0,L)}(z) \, \alpha \nu \Big(\frac{\vx}{\ell} \Big) \Big].
\label{eq:F12}
\end{equation}
It is a random perturbation of the constant  speed $c_o$,  modeled by the random process 
$\nu$. The perturbation extends over the range interval $z \in (0,L)$, as given by the indicator function 
$1_{(0,L)}(z)$. We assume that $\nu(\vecr)$ is a dimensionless stationary random process of
dimensionless argument  $\vecr \in \mathbb{R}^{d+1}$, with zero mean
$\EE \big[ \nu(\vecr) \big] = 0$
and autocovariance
\begin{equation*}
 \EE\big[\nu(\vecr)\nu(\vecr')\big] = \cR(\vecr-\vecr') ,\quad 
\quad \forall \, \vecr, \vecr' \in \mathbb{R}^{d+1}.
\end{equation*}
Moreover, $\nu$ is bounded  
and $\cR$ is integrable, with Fourier transform, the power spectral density
\begin{align}
\widetilde \cR({\vec{\itbf q}}) = \int_{\RR^{d+1}} {\rm d} \vecr\,  \cR(\vecr) e^{-i {\vec{\itbf q}} \cdot \vecr} ,
\label{eq:DefFT}
\end{align}
that is either compactly supported in a ball of radius $O(1)$ in $\mathbb{R}^{d+1}$, or is negligible
outside this ball. The autocovariance is 
normalized by
\begin{equation*}
\int_{\mathbb{R}^{d+1}} {\rm d} \vecr \, \cR(\vecr) = O(1), \qquad \cR({\bf 0}) = O(1).
\end{equation*}
Then, the length scale $\ell$ is the correlation length and the
positive and small dimensionless parameter $\alpha$ quantifies the
typical amplitude (standard deviation) of the fluctuations.

The problem is to characterize the wave field $u(\vx)$ in the scaling regime 
\begin{equation}
\label{eq:sep}
\la   < \ell \sim X \ll L, \qquad \alpha \sim \big({\la}/{L}\big)^{\frac{1}{2}} \ll 1.
\end{equation}
Here $\la = 2 \pi c_o/\om$ is the wavelength and the reference length scale is $L$, which is of the order of the distance of propagation.  We are particularly interested in the coherent (mean) field $\EE[u(\vx)]$ and  the intensity resolved over directions of propagation, the mean Wigner transform of $u(\vx)$. Its evolution in $z$ is governed by the one-way radiative transfer equation that we derive.

\section{Main results}
\label{sect:MR}
Because the interaction of the waves with the random medium depends on
the direction of propagation, we decompose $u(\vx)$ over plane waves, using the Fourier transform with respect to the
transverse coordinates $\bx \in \mathbb{R}^d$ of $\vx = (\bx,z)$,
\begin{equation}
\label{eq:FTr}
\what u(\bka,z) = \int_{\mathbb{R}^d} {\rm d} \bx \, u(\bx,z) e^{-
    i k \bka \cdot \bx}.
\end{equation}
Here $\bka \in \mathbb{R}^d$ is the normalized transverse wave vector, and we suppressed the wavenumber  $k = \omega/c_o$ in the argument
of $\what u$. We show in section \ref{sect:deriv}
that in the scaling regime (\ref{eq:sep}), the field $\what u(\bka,z)$ consists of forward propagating 
waves with longitudinal wavenumber $k \beta(\bka)$, where 
\begin{equation}
\beta(\bka) = \sqrt{ 1-|\bka|^2 } , \quad |\bka| < 1.
\label{eq:WD3}
\end{equation}
The amplitudes of these waves (modes) are denoted by  $a(\bka,z)$. They are complex-valued  $z$-dependent random 
fields which model  wave scattering in the  random medium. 

The wave field $u(\vx)$ is given by the Fourier synthesis of the modes, the plane waves with  wave vector
$k\vec{\bka} = k( \bka , \beta(\bka) )$, 
 \begin{align}
 u(\vx) =\int_{|\bka|<1} \frac{{\rm d}(k \bka)}{(2 \pi)^d}  
  \frac{a(\bka,z)}{\beta^{\frac{1}{2}}(\bka) }   e^{i k 
    \vec{\bka} \cdot \vx}  , \quad \vx = (\bx,z) ,
   \label{eq:synthesis}
    \end{align}
where have used the notation ${\rm d}(k \bka) = k^d {\rm d} \bka$ for the infinitesimal volume in $\mathbb{R}^d$.
The mode amplitudes are normalized by the factors $\beta^{\frac{1}{2}}(\bka)$ in order to simplify the formulae that follow
\footnote{In particular these factors ensure that the energy fluxes of the plane wave modes through the planes $z=$ constant are $|a(\bka,z)|^2$.}.
In the scaling regime (\ref{eq:sep}) the mode amplitudes form a Markov process whose statistical moments can be characterized explicitly, as explained  in subsection \ref{sect:DL}.
 Here we describe the expectation of $a(\bka,z)$, which defines the coherent field, and 
its second moments, which define the mean Wigner transform of $u(\vx)$.

The mean mode amplitudes are 
\begin{equation*}
\EE[a(\bka,z)] = a_o(\bka) \exp[{\rm Q}(\bka) z], 
\end{equation*}
where $a_o(\bka)$ are the  amplitudes in the homogeneous medium, defined in Eq.~(\ref{eq:WD8bp}) by the source excitation. 
The effect of the random medium 
is in the complex exponent
\begin{align}
\nonumber 
{\rm Q}(\bka) &=-\frac{k^2 \alpha^2 \ell^{d+1}}{4 } \int_{|\bka'|<
   1}  \frac{{\rm d} (k \bka')}{(2 \pi)^d}\, \frac{1}{\beta(\bka)\beta(\bka')} 
\\
 &  \times 
 \int_0^\infty  
  {\rm d} \zeta \int_{\mathbb{R}^d} {\rm d} \br \,  \cR (\br, \zeta ) e^{- 
   i k \ell \big(\bka-\bka',\beta(\bka)-\beta(\bka')\big)\cdot (\br,\zeta)} .
\label{eq:C5}
\end{align}
Since $\cR$ is even, the real part of ${\rm Q}(\bka)$ is determined by the power spectral density $\widetilde{\cR}$
defined in \eqref{eq:DefFT},
which is non-negative by Bochner's theorem \footnote{ Bochner's theorem states that a function
is an autocovariance function of a stationary process if and only if its Fourier transform is nonnegative \cite{gihman}.}.
Thus  ${\rm Re}\big[{\rm Q}(\bka)\big] < 0$,
and the mean amplitudes decay exponentially in $z$,
with the  decay rate  
\begin{align}
\frac{1}{{\mathcal S}(\bka)} =  -  \mbox{Re}\big[ {\rm Q}(\bka) \big].
\label{eq:C5b}
\end{align}
The length ${\mathcal S}(\bka) $
 is the scattering mean free path.
By choosing the magnitude $\alpha$ of the fluctuations  as in (\ref{eq:sep}), we have $L \sim \mathcal{S}(\bka)$, so the decay with $z$ is significant in our scaling regime.  It is the manifestation of the randomization of the wave, due to scattering in the  medium. 

The strength of the random  fluctuations of the mode amplitudes is described by the energy density (Wigner transform) 
\begin{align}
{\mathcal W}(\bka,\bx,z) = \int \frac{{\rm d}(k {\itbf q})}{(2 \pi)^d}
\, \exp\Big[i k \itbf q \cdot \big(\nabla \beta(\bka) z + \bx\big)\Big]  ~~\nonumber \\  \times 
  \EE \left[   a \big(\bka
    + \frac{{\itbf q}}{2},z\big) \overline{    a\big(\bka-\frac{ {\itbf q}}{2},z\big)}\right],
 \label{eq:T6:0}
\end{align}
where the bar denotes complex conjugate and the integral is over all $\bq \in \RR^d$ such that $|\bka \pm \bq/2|<1$. 
 The Wigner transform satisfies the transport equation 
\begin{align}
\nonumber & \partial_z {\mathcal W} (\bka,\bx,z)  - \nabla \beta(\bka) \cdot \nabla_\bx
  {\mathcal W}  (\bka,\bx,z) =\\
  & \int_{|\bka'| < 1} \frac{{\rm d} (k \bka')}{(2 \pi)^d}\,
{\mathcal Q}(\bka,\bka')  \big[  {\mathcal W}(\bka',\bx,z)   - {\mathcal W}(\bka,\bx,z)   \big] ,
\label{eq:T7:0}
\end{align}
for $z>0$, with differential scattering cross section 
\begin{align}
\hspace{-0.09in} {\mathcal Q}(\bka,\bka') = \frac{k^2 \alpha^2 \ell^{d+1} }{4 \beta(\bka)\beta(\bka') } \widetilde  \cR 
\Big( k \ell \big(\bka-\bka'\big),k \ell \big(\beta(\bka)-\beta(\bka')\big)\Big).
  \label{eq:DCS}
\end{align}
The total scattering cross section is 
\begin{align}
\Sigma(\bka) = \int_{|\bka'| < 1} \frac{{\rm d} (k \bka')}{(2 \pi)^d}\,
{\mathcal Q}(\bka,\bka') = - 2 \mbox{Re} \big[{\rm Q}(\bka) \big] = \frac{2}{{\mathcal S}(\bka)}.
\end{align}
Equation (\ref{eq:T7:0}) looks like the radiative transfer equation, except that it is an initial value problem in $z$,  with ${\mathcal W}(\bka,\bx,z=0)$ given 
by the Wigner transform of mode amplitudes $a_o(\bka)$ in the homogeneous medium.
As we show in subsection \ref{sect:ResA} it is in fact a general form of the one-way  radiative transfer equation
introduced recently in the biomedical imaging literature \cite{kim15}.  We also establish in subsection 
\ref{ap:RT} the connection between equation (\ref{eq:T7:0}) and the standard radiative transfer theory: 
We show that Eq.~(\ref{eq:T7:0}) can be obtained heuristically from the standard radiative transfer equation by
applying a forward scattering approximation.  Such a calculation is heuristic, because the standard radiative transfer equation has no rigorous derivation \cite{schotland15}, whereas 
Eq. (\ref{eq:T7:0}) is derived here from first principles. 
The connection to the It\^{o}-Schr\"{o}dinger model is  
in subsection \ref{app:C}:
We show that Eq.~(\ref{eq:T7:0}) can be reduced to the transport equation in the paraxial
geometry by taking the limit of very small angles.
Therefore Eq.~(\ref{eq:T7:0}) can be seen as a bridge between the radiative transfer and paraxial approximation 
regimes. 

\subsection{Connection with the one-way radiative transfer equation}
\label{sect:ResA}
The one-way radiative transfer equation was proposed recently in \cite{kim15} for the application 
of diffusion optical tomography in forward-peaked scattering media. The equation is stated in \cite{kim15} in two dimensions $(d+1 = 2)$, 
\begin{align}
\label{eq:owrte}
\sin \theta \partial_zI +\cos \theta \partial_x I= \mu_s \int_0^\pi p(\theta-\theta') \big[ I(\theta')-I(\theta) \big]  {\rm d} \theta' ,
\end{align}
for $I(\theta,x,z)$ the light intensity at position $(x,z)$ in the direction $(\cos \theta,\sin \theta)$, with $\theta \in [0,\pi]$. 
The coefficient   $\mu_s$  is the total scattering cross section and the scattering phase function $p(\theta-\theta')$ is chosen of the Henyey-Greenstein form \cite{heino03,kim15},
\begin{align}
\label{eq:hg}
p(\theta-\theta') = \frac{1}{2\pi} \frac{1-g^2}{1+g^2 -2g \cos (\theta-\theta')}, 
\end{align}
satisfying $\int_0^{2\pi} p(\theta) d\theta=1$. Parameter $g\in (0,1)$ is the anisotropy factor and it is argued that the one-way radiative transfer equation is valid when $g \sim 1$, so scattering is forward-peaked.

The light intensity $I$ is in fact the Wigner transform ${\mathcal W}$
introduced in (\ref{eq:T6:0}), with $\bka=\cos \theta \in (-1,1)$.
Indeed, in statistically isotropic media, i.e., $\cR(\vec\bx)=\cR_{\rm iso}(|\vec\bx|)$,
we obtain from (\ref{eq:T7:0}) (multiplied by $\sin \theta$),  using that $\beta(\bka)=\sin \theta$ and $\nabla \beta(\bka)= - \cot \theta $, 
\begin{align}
\nonumber
& \sin \theta \partial_z {\mathcal W} +\cos \theta  \partial_x
  {\mathcal W}= 
  \frac{k^3 \ell^2 \alpha^2 }{4 } \\
 & \times
  \int_0^\pi {\rm d} \theta'  \,  \breve \cR_{\rm iso}\Big(k \ell \sqrt{2 (1-\cos(\theta-\theta')} \Big)
   \big[  {\mathcal W}(\theta')  - {\mathcal W}(\theta)  \big] ,
  \label{eq:T7b}
\end{align}
with 
\begin{align}
 \breve \cR_{\rm iso}(q) = \int_0^\infty {\rm d}s \, s \cR_{\rm iso}(s) J_0(qs)  .
 \label{eq:defRiso}
\end{align}
This is exactly (\ref{eq:owrte}) with the identification:
\begin{equation}
\mu_s p(\theta-\theta') =   \frac{k^3 \ell^2 \alpha^2 }{4 }   \breve \cR_{\rm iso} \Big(k \ell \sqrt{2 (1-\cos(\theta-\theta')} \Big) .
\label{eq:identif}
\end{equation}
The scattering phase function (\ref{eq:hg}) is a particular case of \eqref{eq:identif}, corresponding to a Lorentzian
for  $ \breve \cR_{\rm iso}$, that is 
\begin{equation}
\label{eq:identif1}
\breve  \cR_{\rm iso} (q) = \frac{\breve \cR_o}{1+q^2}.
\end{equation}
This corresponds (through (\ref{eq:defRiso}) and \cite[formula 6.521.2]{grad}) 
to an autocovariance function of the form $\cR_{\rm iso}(s) =\breve \cR_o K_0(s)$,
where $K_0$ is the Bessel function of the second kind of order zero.
This is the zeroth von K\'arm\'an correlation function \cite{klimes}.
It has a logarithmic divergence at $s=0$, 
which can be regularized by introducing an ultraviolet cutoff in (\ref{eq:identif1}).
By substituting (\ref{eq:hg}) and (\ref{eq:identif1}) into (\ref{eq:identif}) we obtain 
the anisotropy parameter and total scattering cross section 
\begin{align*}
g= 1+\frac{1}{2(k\ell)^2} - \frac{1}{k \ell} \sqrt{1+\frac{1}{4(k \ell)^2}} , \quad ~
\mu_s=  \Big(\frac{1-g}{1+g}\Big)  \frac{\pi k^3 \ell^2 \alpha^2 \breve \cR_o}{2}   .
\end{align*}
The validity condition $g\sim1$  in \cite{kim15} is equivalent
to  $\lambda < \ell$. 
This completes the proof that (\ref{eq:owrte}) is a special case of 
our Eq.~(\ref{eq:T7:0}). 
It justifies the model (\ref{eq:owrte}), as our results in this paper show that it can be rigorously derived from the wave equation in random media,
in the scaling regime (\ref{eq:sep}).

\subsection{Connection to the radiative transfer theory}
\label{ap:RT}
To connect our transport equation (\ref{eq:T7:0}) to the  standard radiative transfer theory in   random media 
\cite{chandra,ryzhik96,bal00}, we let $d+1 = 3$ and adhere to the notation in
\cite{ryzhik96}. 
Following  \cite[Eq.~(3.42)]{ryzhik96}, we define 
\begin{equation*}
f (\vk,\vx) =  {\pi}  \Big[ - \frac{i }{k} \frac{\vk}{|\vk|}  \cdot   {\vec \nabla}_{\vx} u(\vx)+ u(\vx)  \Big], 
\end{equation*}
where we use a different constant of proportionality than in
\cite{ryzhik96}, to simplify the relation in \eqref{eq:RT6}.  The Wigner
transform $ W(\vk,\vx)$ introduced in \cite[Eq.~(3.41)]{ryzhik96}  is 
\begin{equation}
W (\vk,\vx) = \int_{\RR^3} \frac{{\rm d} \vy}{(2 \pi)^3} \, f \Big(\vk, \vx -
\frac{\vy}{2}  \Big) \overline{f\Big( \vk ,\vx + \frac{
    \vy}{2} \Big)} e^{ i \vk \cdot \vy}, 
\label{eq:RT5}
\end{equation}
and satisfies the transport equation  \cite[Eq.~(4.38)]{ryzhik96}
\begin{align}
\nonumber
 \vna_{\vk} \om(\vk) \cdot \vna_{\vx} W(\vk,\vx)
= \int_{\RR^3}  {\rm d} \vk' \sigma( \vk,\vk')  W(\vk',\vx) 
\\-
\Sigma(\vk) W(\vk,\vx),
\label{eq:RTT}
\end{align}
with dispersion relation
$
\om(\vk) = c_o |\vk|,
$
and integral kernel,  the
differential scattering cross-section, 
\begin{align}
\sigma( \vk,\vk')   = \frac{\pi c_o^2 k^2 \ell^3 \alpha^2}{2(2
  \pi)^3} \widetilde \cR \big[ \ell (\vk-\vk' ) \big]
\delta[\om(\vk)-\om(\vk')]   .
\label{eq:RT9}
\end{align}
The scalar $\Sigma(\vk)$ is the total scattering cross section 
\begin{align}
\Sigma(\vk)  = \int_{\RR^3}  {\rm d} \vk'\, \sigma(\vk,\vk')  . 
\label{eq:RT9b}
\end{align}

Substituting (\ref{eq:synthesis}) into (\ref{eq:RT5}),
 we obtain after some
algebraic manipulations that
\begin{equation}
\label{eq:RT6}
W(\vk,\vx) = \frac{\delta \big[ \mk_z - k \beta(\bmk/k)
     \big]}{\beta(\bmk/k)} {\mathcal W}(\bmk/k,\bx,z),
\end{equation}
with ${\mathcal W}$ the Wigner transform   (\ref{eq:T6:0}).
The Dirac factor in Eq.~(\ref{eq:RT6}) expresses the fact that in our scaling regime,
in which the wave field has the form (\ref{eq:synthesis}),
the forward scattering approximation is valid and the intensity resolved over directions of propagation
is supported on the wave vectors $\bmk$ with positive $\mk_z$.
Next we rewrite the three  terms of (\ref{eq:RTT}) to show that the equation   is equivalent to \eqref{eq:T7:0}.

1) Since \eqref{eq:RT6} gives that $W(\vk,\vx)$ is supported at vectors $\vk$ of
the form $\vk = k \vka$, with $\vka = (\bka,\beta(\bka))$, the operator on the left hand side of (\ref{eq:RTT}) is 
\begin{equation*}
{\vna}_{\vk} \om(\vk) \cdot \vna_{\vx}   = c_o \beta(\bka) \big[ \partial_z - \nabla
  \beta(\bka) \cdot \nabla_{\bx} \big],
\end{equation*}
and   we obtain that 
\begin{align}
\nonumber
{\vna}_{\vk}  \om(\vk) \cdot \vna_{\vx} W(\vk,\vx) =c_o \delta \big[\mk_z
   - k \beta(\bmk/k)\big]
 \\
 \times \big[ \partial_z - \nabla \beta(\bmk/k) \cdot
   \nabla_{\bx} \big] {\mathcal W}(\bmk/k,\bx,z). \label{eq:RT11}
\end{align}

2) The integral kernel in  (\ref{eq:RTT}) 
is supported at  $\vk' = k \vka'$, with $\vka' =
(\bka',\beta(\bka'))$, by \eqref{eq:RT6}, so the Dirac distribution in (\ref{eq:RT9}) is
\begin{align*}
\delta \big[ \om(\vk)-\om(k \vka')\big] &
= \frac{\delta
  \big[ \mk_z - k \beta(\bmk/k)\big]}{c_o \beta(\bmk/k)}.
\end{align*}
Thus,  we have
\begin{align}
&\hspace*{-0.05in}
\int_{\RR^3}  {\rm d} \vk' \sigma( \vk,\vk')  W(\vk',\vx) = \frac{c_o k^2 \ell^3\alpha^2}{4 } \delta \big[\mk_z
  - k \beta(\bmk/k)\big] \nonumber \\ 
  &\hspace*{-0.05in}
\times \hspace*{-0.05in} \int_{|\bka'|< 1} \hspace{-0.05in}\frac{{\rm d} (k
  \bka')}{(2 \pi)^2}\, \frac{\widetilde \cR \big[ \ell(\bmk-k\bka'), k\ell(\beta(\bmk/k)-\beta(\bka'))
\big]  }{\beta(\bmk/k)\beta(\bka')}{\mathcal W}(\bka', \bx,z)
   ,\label{eq:RT12}
\end{align}
where $|\bka'| < 1$ because we have only propagating waves.

3) From (\ref{eq:RT9b}) we find that
\begin{align*}
\Sigma(\vk)  = \frac{ c_o^2 k^2 \ell^3 \alpha^2}{4(2
  \pi)^2}  \int_{\RR^3}  {\rm d} \vk'  \delta \big[ \omega(\vk')-\omega(\vk)\big]
  \widetilde \cR\big( \ell(\vk-\vk')\big)  ,
\end{align*}
so for $\vk=k(\bka,\beta(\bka))$,
\begin{align}
\nonumber
&\Sigma(\vk) W(\vk,\vx) = \frac{ c_o k^2 \ell^3 \alpha^2}{4}  \delta \big[ \mk_z - k \beta(\bmk/k)
     \big]  \int_{|\bka'| < 1} \frac{{\rm d} (k \bka')}{(2
  \pi)^2}  \\
 &  \times \frac{ \widetilde \cR\big[ \ell(\bmk-k\bka'),\ell k(\beta(\bmk/k)-\beta(\bka'))
\big]}{\beta(\bmk/k)\beta(\bka')}
  {\mathcal W}(\bmk/k,\bx,z) .
\label{eq:RT12b}
\end{align}

Finally, substituting (\ref{eq:RT11}), (\ref{eq:RT12}), and (\ref{eq:RT12b}) into the transport equation (\ref{eq:RTT}) satisfied by
$W$,  we obtain that the Wigner transform ${\mathcal W}$ satisfies the transport equation (\ref{eq:T7:0}).
This completes the proof that Eq.~(\ref{eq:T7:0}) can be obtained from the standard radiative transfer equation (\ref{eq:RTT})
by applying a forward scattering approximation. However, as stated before, there is no rigorous derivation of the standard  radiative transfer equation from the wave equation in random media.  In this paper we obtain a rigorous derivation of Eq.~(\ref{eq:T7:0}) from the wave equation in random media, in the scaling regime (\ref{eq:sep}).

\subsection{Connection to the paraxial theory}
\label{app:C}
It is shown in \cite{garnier2009paraxial} that if $\la \ll \ell \ll L$ so that the medium Fresnel number ${\ell^2}/({\la L}) \sim 1$,  
and if the standard deviation $\alpha$ of the fluctuations is small so that $\alpha^2 \sim {\la^2}/({\ell L})$, 
then the inverse Fourier transform of the mode amplitudes, denoted by $a_{\rm pa}(\bka,z)$,
\begin{align*}
\check{a}_{\rm pa}(\bx,z) = \int_{|\bka|<1} \frac{{\rm d} (k \bka)}{(2\pi)^d} \, a_{\rm pa}(\bka,z) e^{i k \bka \cdot \bx} ,
\end{align*}
satisfies the random paraxial wave equation 
 (or It\^{o}-Schr\"{o}dinger model) \cite{garnier2009paraxial}
\begin{align}
\hspace{-0.09in}{\rm d} \check{a}_{\rm pa}(\bx,z) = \frac{i}{2k} \Delta_\bx \check{a}_{\rm pa}(\bx,z){\rm d}z +\frac{ik}{2} \check{a}_{\rm pa}(\bx,z) \circ {\rm d}B(\bx,z).
\label{eq:ItoSchr}
\end{align}
Here $B$ is the Brownian field i.e., a Gaussian process with mean zero and covariance
\begin{align*}
\EE [ B(\bx,z)B(\bx',z') ]=  \alpha^2 \ell \min(z,z') {\mathscr C}\Big(\frac{\bx-\bx'}{\ell} \Big), \quad 
{\mathscr C}(\br) = \int_{-\infty}^\infty {\rm d} \zeta \,  \cR(\br,\zeta ).
\end{align*}
The symbol  $\circ$ stands for the Stratonovich integral.  
This integral is the suitable form  of stochastic integral for the It\^{o}-Schr\"{o}dinger model as shown in 
\cite{garnier2009paraxial}, and as could be predicted by the general Wong-Zakai theorem \cite{wongzakai}. 
Alternatively, we can   characterize $a_{\rm pa}(\bka,z)$ as the solution of 
\begin{align*}
\hspace{-0.09in}{\rm d} \check{a}_{\rm pa}(\bx,z) = \frac{i}{2k} \Delta_\bx \check{a}_{\rm pa}(\bx,z){\rm d}z +\frac{ik}{2} \check{a}_{\rm pa}(\bx,z)   {\rm d}B(\bx,z) - \frac{k^2 \ell \alpha^2   {\mathscr C}({\bf 0}) }{8}\check{a}_{\rm pa}(\bx,z) {\rm d} z,
\end{align*}
where the stochastic integral is now understood in the usual It\^o's form.

The derivation of \eqref{eq:ItoSchr} from the wave equation in random media,
given in  \cite{garnier2009paraxial}, involves two main steps: 
first show that the forward scattering approximation is valid; second show that the effect of the 
fluctuations of the random medium on the wave field can be captured in distribution by a white noise (in $z$) 
model.

Using the  It\^o-Schr\"odinger model  \eqref{eq:ItoSchr} we find by It\^o's formula that  the mean field $\check{A}_{\rm pa}(\bx,z) =  \EE [\check{a}_{\rm pa}(\bx,z)]$
satisfies 
\begin{align*}
\partial_z  \check{A}_{\rm pa}(\bx,z) = \frac{i}{2k} \Delta_\bx \check{A}_{\rm pa}(\bx,z)- \frac{k^2 \ell\alpha^2 {\mathscr C}({\bf 0})}{8} \check{A}_{\rm pa}(\bx,z)  .
\end{align*}
It decays with $z$ on the scale 
\begin{align*}
{\mathcal S}_{\rm pa} =  \frac{8}{k^2 \ell \alpha^2 {\mathscr C}({\bf 0})} = \frac{8}{k^2 \ell \alpha^2   \int_{-\infty}^\infty {\rm d} \zeta \,  \cR({\bf 0},\zeta)},
\end{align*}
which corresponds to the scattering mean free path $\cS(\bka)$ defined by (\ref{eq:C5}-\ref{eq:C5b}), for $\la \ll \ell$ and 
$|\bka| = O({\la}/{\ell})$.

The Wigner transform is 
\begin{align*}
\nonumber
W_{\rm pa} (\bmk,\bx,z) = \int_{\RR^d} {\rm d}\by \, e^{ i \bmk \cdot \by } \EE \Big[
 \check{a}_{\rm pa}\big(\bx-\frac{\by}{2},z\big) \overline{\check{a}_{\rm pa}\big(\bx+\frac{\by}{2} ,z\big)} \Big]  \\
 = \int_{\RR^d}\frac{{\rm d}(k\bq)}{(2\pi)^d} \, e^{ i k\bq\cdot \bx } \EE \Big[
 a_{\rm pa}\big( \frac{\bmk}{k}+\frac{\bq}{2},z\big) \overline{a_{\rm pa}\big( \frac{\bmk}{ k}-\frac{\bq}{2},z \big)} \Big] ,
\end{align*}
which corresponds to (\ref{eq:T6:0}) for $\bmk=k\bka$ and $|\bka| = O({\la}/{\ell})$.
Using It\^o's formula it is shown in \cite{garnier2009paraxial} to satisfy the transport equation 
\begin{align}
\nonumber
&\hspace*{-0.075in} \partial_z W_{\rm pa} +\frac{1}{k} \bmk \cdot \nabla_\bx W_{\rm pa} =\\ 
&\hspace*{-0.075in} \frac{k^2 \ell^{d+1} \alpha^2}{4} \int_{\RR^d} \frac{{\rm d} \bmk'}{(2\pi)^d}\, \widetilde \cR \big(\ell (\bmk-\bmk') ,0 \big)
\big[ W_{\rm pa}(\bmk') - W_{\rm pa}(\bmk) \big] ,
\label{eq:rtepar}
\end{align}
with differential scattering cross section  
\begin{align*}
{\mathcal Q}_{\rm pa}(\bmk,\bmk') &= \frac{k^2 \ell^{d+1} \alpha^2 }{4} \widetilde \cR \big(\ell (\bmk -\bmk') ,0 \big)
\end{align*}
corresponding to (\ref{eq:DCS}) for $\bmk=k\bka$, $\bmk'=k\bka'$, and $|\bka|, |\bka'| = O({\la}/{\ell})$.

This establishes the connection between Eq. (\ref{eq:T7:0}) and the transport equation 
(\ref{eq:rtepar}) derived in \cite{garnier2009paraxial}. Together with the result in section \ref{ap:RT} it completes the proof that Eq. (\ref{eq:T7:0}) is a bridge between the radiative transfer and paraxial approximation regimes.

We end the section with the note that, as shown for instance in \cite[Chapter 13]{ishimaru},
the radiative transfer equation in the white-noise paraxial regime (\ref{eq:rtepar})
can also be derived heuristically from the standard radiative transfer equation in the ``approximation of large particles", 
or equivalently in the ``small angle approximation", which corresponds to a random medium
with large correlation radius.

\section{Derivation of results}
\label{sect:deriv}
To derive the transport equation (\ref{eq:T7:0}) from the wave equation, we use multiscale analysis and 
probabilistic limit theorems. The asymptotic regime of separation of scales (\ref{eq:sep}) is defined in terms of 
three small dimensionless parameters 
\begin{equation}
\eps= \frac{\la}{L}, \qquad \gamma = \frac{\la}{\ell} ,\qquad {\eta} = \frac{\la}{X},
\label{eq:F17}
\end{equation}
ordered as
\begin{equation*}
0 < \eps \ll \gamma \sim {\eta} < 1,
\end{equation*}
and the standard deviation $\alpha $ of the fluctuations of the random medium is of order $\eps^{\frac{1}{2}}$.  We begin 
with the wave decomposition, and obtain a stochastic system of differential equations satisfied by the mode amplitudes.
We consider both forward and backward going waves, but then show that we can neglect the backward waves
in the limit $\eps \to 0$ (subsection \ref{sect:DL}).  The $\eps \to 0 $ limit of the mode amplitudes defines the Markov process whose expectation and Wigner transform are described in section \ref{sect:MR}. 

\subsection{Scaled equation}
\label{sect:multi}%
We let ${L}$ be the reference length scale, which is similar to the distance of propagation, and 
introduce the scaled length variables
$
\bx' = {\bx}/(\eps L), $ $z' = {z}/{L},$ $ L' = {L}/{L} = 1$, 
$\ell' ={\ell}/{L} = {\eps}/{\gamma}$ and $
X' = {X}/{L} = {\eps}/{\eta}$.
The scaled standard deviation is $\alpha'=\alpha / \eps^{1/2}$.
The scaled wavenumber is $k' = {k L}{\eps}=2\pi$.

Let us denote the wave field by $u^\eps$. Substituting in
(\ref{eq:F0}) and dropping all the primes, as all
the variables are scaled henceforth, we obtain
\begin{align}
\bigg\{\partial_z^2  + \frac{1}{\eps^2} \Delta_\bx + \frac{k^2}{\eps^2}
 \Big[1+  \eps^{\frac{1}{2}} \alpha \nu \Big(\gamma \bx, \frac{\gamma z}{\eps}\Big)\Big]
 \hspace*{-0.02in}
 \bigg\}
 u^\eps(\bx,z)
=
-  \frac{1}{\eps}F\big( {\eta} {\bx} \big)
\delta(z), 
\label{eq:F19}
\end{align}
for $0 \le z \le L$. At ranges
$z < 0$ and $z > L$ the equations are simpler, as the term 
involving the process $\nu$ vanishes.   
Since the wave field depends  linearly on  the source,  we scaled $F$ 
by $1/\eps$ to obtain an order one result in the  limit $\eps \to 0$. 

\subsection{Wave decomposition}
\label{sect:WD}
We decompose the field $u^\eps(\bx,z)$ in plane waves using the  Fourier transform  with respect to $\bx \in \mathbb{R}^d$, as in \eqref{eq:FTr}:
\begin{equation}
\what u^\eps(\bka,z) = \int_{\mathbb{R}^d} {\rm d} \bx \, u^\eps(\bx,z) e^{-
    i k \bka \cdot \bx}.
\end{equation}
The transformed field $\what u^\eps(\bka,z)$
is a superposition of  forward and backward going waves (modes) along $z$, as explained next. To ease the explanation 
we begin with the reference case in the homogeneous medium, and then consider the random medium.

\subsubsection{Homogeneous media}
\label{sect:WDH}
The transformed field in homogeneous media $\what
u_o^\eps(\bka,z)$ satisfies the ordinary differential equation 
\begin{align}
\partial_z^2 \what u_o^\eps(\bka,z) + \frac{k^2}{\ep^2} \beta(\bka)^2\what u_o^\eps(\bka,z)
=
- \frac{1}{\eps {\eta}^d} \what F\Big( \frac{k\bka}{{\eta}} \Big)
\delta(z), \label{eq:F19b}
\end{align}
with $\beta(\bka)$ defined in (\ref{eq:WD3}) and $\what F$ the Fourier transform of $F$,
\begin{align}
\what F({\itbf q}) = \int_{\RR^d} F(\br) e^{- i {\itbf q} \cdot \br} {\rm d} \br .
\label{def:fourierF}
\end{align}
The solution is  outgoing and bounded away from the source, and it is given explicitly, for $z \ne 0$, by 
\begin{align}
 \label{eq:WD4}
\hspace{-0.1in}\what u_o^\ep(\bka,z) &= 
\frac{a_o (\bka )}{\beta^{\frac{1}{2}}(\bka)}
e^{\frac{ik}{\ep} \beta(\bka)z}  1_{(0,\infty)} (z) +
\frac{b_o(\bka)}{\beta^{\frac{1}{2}}(\bka)}
e^{-\frac{ik}{\ep} \beta(\bka)z} 1_{(-\infty,0)}(z)  .
\end{align}
Thus, the wave field 
\[
u_o^\eps(\bx,z) = \int_{|\bka| < 1}  \frac{{\rm d} (k \bka)}{(2 \pi)^d} \what u_o^\eps(\bka,z) e^{i k \bka \cdot \bx}
\]
is a synthesis of plane waves with wave vectors 
$k\big(\bka,\pm \beta(\bka)\big)$.  The plus sign corresponds to forward going waves, and the negative sign to backward going
waves. The amplitudes are determined by the
jump conditions at the source
\begin{align*}
\what u_o^\eps(\bka, 0+)-\what u_o^\eps(\bka,0-) &= 0 , \\
\partial_z \what u_o^\ep(\bka, 0+)-\partial_z \what u_o^\eps(\bka,0-) &=
-\frac{1}{\eps {\eta}^d} \what F\Big( \frac{k\bka}{{\eta}}\Big),
\end{align*}
which gives
\begin{align}
a_o(\bka) = b_o(\bka) = 
\frac{i }{2 k {\eta}^d \beta^{\frac{1}{2}}(\bka) } \what F
  \Big( \frac{k\bka}{{\eta}}\Big) .
   \label{eq:WD8bp}
\end{align}

The radius of the support of $\what F (\bq)$ is one, so the scaling parameter $\eta$ 
controls  the support in $\bka$ of the wave modes generated  by the source i.e., the opening 
angle of the initial wave beam. 
Consistent with  (\ref{eq:sep}) and \eqref{eq:F17}, we assume henceforth that  
\begin{equation}
\frac{\eta }{k} <1.
\label{eq:ASG1}
\end{equation}
so that in  (\ref{eq:WD8bp}) we have $ |\bka| \leq {\eta} / k < 1 $. Then
$\beta(\bka)$ defined by (\ref{eq:WD3}) is real valued, and there are
no evanescent waves in the decomposition (\ref{eq:WD4}).

\subsubsection{Random media}
\label{sect:WDR}
The field $\what u^\eps(\bka,z)$ in the random medium satisfies
the equation
\begin{align}
\partial_z^2 \what u^\eps+ \frac{k^2}{\eps^2} \beta(\bka)^2\what u^\eps
   +  1_{(0,L)}(z)
\,   \mathcal{M}^\eps\what u^\ep  
= - \frac{1}{\eps {\eta}^d} \what F\Big( \frac{k \bka}{{\eta}} \Big)
\delta(z),
\label{eq:WD10}
\end{align}
derived from (\ref{eq:F19}), with radiation conditions
at $z < 0$ and $z > L$, and source conditions at $z = 0$. The leading
$O(1/\eps^2)$ term in the right hand side is the same as in the homogeneous medium, so 
we can use a similar wave decomposition to that in section \ref{sect:WDH}. 
The random perturbation  is in the operator $\mathcal{M}^\eps$
defined by
\begin{align*}
\mathcal{M}^\eps
      \what u^\eps(\bka,z) =&
  \frac{ik \alpha}{\eps^{\frac{1}{2}}\gamma^d} \int \frac{{\rm d} (k
    \bka')}{(2 \pi)^d} \, \frac{\what \nu \Big(\frac{k(\bka
    -\bka')}{\gamma}, \frac{\gamma z}{\eps}\Big)}{[\beta(\bka)\beta(\bka')]^{\frac{1}{2}}} 
 \what u^\ep(\bka',z) ,
\end{align*}
where $\what \nu$ is the Fourier transform 
of $\nu$  with respect to the first argument in $\mathbb{R}^d$ as in (\ref{def:fourierF}).

The wave decomposition is 
\begin{align*}
\hspace*{-0.085in}
a^\eps(\bka,z) &= \frac{1}{2} 
\Big( \beta(\bka)^{\frac{1}{2}} \what u^\ep(\bka,z)  +\frac{\eps}{ik \beta(\bka)^{\frac{1}{2}}} \partial_z \what u^\ep(\bka,z)
\Big)  e^{- \frac{i k}{\eps} \beta(\bka) z}  ,\\
\hspace*{-0.085in}
b^\eps(\bka,z) &= \frac{1}{2} 
\Big( \beta(\bka)^{\frac{1}{2}} \what u^\ep(\bka,z)  -\frac{\eps}{ik \beta(\bka)^{\frac{1}{2}}} \partial_z \what u^\ep(\bka,z)
\Big)  e^{ \frac{i k}{\eps} \beta(\bka) z}  ,
\end{align*}
so that we can write as in the homogeneous medium
\begin{align}
\label{eq:WD14}
\what u^\eps(\bka,z)  =& \frac{1}{\beta(\bka)^{\frac{1}{2}}}
\Big(  a^\eps(\bka,z)
 e^{\frac{i k}{\eps} \beta(\bka) z}  +
  b^\eps(\bka,z) e^{-\frac{i k}{\ep}
  \beta(\bka) z} \Big) , \\
\partial_z \what u^\eps(\bka,z)  =& \frac{ i k \beta(\bka)^{\frac{1}{2}}}{\eps}
\Big(  a^\ep(\bka,z)
 e^{\frac{i k}{\eps} \beta(\bka) z}  -
  b^\eps(\bka,z) e^{-\frac{i k}{\eps}
  \beta(\bka) z} \Big) .  
\label{eq:WD14b}
\end{align}
The forward and backward going wave amplitudes $a^\ep(\bka,z)$ and $b^\ep(\bka,z)$
are no longer constant, but random fields due to scattering in the range interval $z \in (0,L)$. 
The  medium  is homogeneous outside this interval and we have the radiation
 conditions 
\begin{align}
  a^\eps(\bka,z) = 0  \mbox{ if } z <
  0 ~~\mbox{ and }~~ b^\eps(\bka,z) = 0 
  \mbox{ if } z\ge L . 
  \label{eq:WD15}
\end{align}
Moreover, $ a^\eps(\bka,z) = a^\eps(\bka,L) $  for $z > L$,
and
$  b^\eps(\bka,z) = b^\eps(\bka,0-)$ for $z < 0$.

The jump conditions at the source are as in section \ref{sect:WDH}, and give
\begin{align}
  a^\eps(\bka,0+) = a_o(\bka) 
~~\mbox{ and }~~
  b^\eps(\bka,0-) = b_o(\bka) + b^\eps(\bka,0+). \label{eq:WD17}
\end{align}
As expected,  the forward going waves leaving the
source are the same as in the homogeneous medium, because the scattering effects 
in the random medium manifest only at long distances of propagation.  The waves at $z < 0$ are given by the
superposition of those emitted by the source, modeled by $b_o(\bka)$, 
and the waves backscattered by the random
medium, modeled by $b^\ep(\bka,0+)$. 

To determine the amplitudes in the random medium, we substitute
equations (\ref{eq:WD14})-(\ref{eq:WD14b}) into (\ref{eq:WD10}). We obtain that 
\begin{align}
\nonumber
\partial_z \begin{pmatrix}
{a}^\eps(\bka,z) \\
{b}^\eps(\bka,z) 
\end{pmatrix}=
\frac{i k \alpha}{2 \gamma^d \eps^{\frac{1}{2}}} \int \frac{{\rm d} (k \bka')}{(2\pi)^d}\, 
\widehat{\nu} \Big( \frac{k(\bka-\bka')}{\gamma},
\frac{\gamma z}{\ep}\Big)\\
\times \boldsymbol{\Gamma}\Big(\bka,\bka',\frac{z}{\eps} \Big)  \begin{pmatrix}
{a}^\ep(\bka',z) \\
{b}^\ep(\bka',z) 
\end{pmatrix} ,
\label{eq:twopointbvp}
\end{align}
in $z \in (0,L)$, with boundary conditions (\ref{eq:WD15})-(\ref{eq:WD17}). 
We are interested in the propagating waves,
corresponding to $|\bka| < 1$ in (\ref{eq:twopointbvp}), and we
explain in section \ref{sect:DL} that in our regime the evanescent
waves may be neglected. The
$2 \times 2$ complex matrices 
\begin{equation}
\boldsymbol{\Gamma} (\bka,\bka',\zeta ) =
\begin{pmatrix}
\Gamma^{aa} (\bka,\bka',\zeta )  &\Gamma^{ab} (\bka,\bka',\zeta ) \\
\Gamma^{ba} (\bka,\bka',\zeta )  &\Gamma^{bb} (\bka,\bka',\zeta ) 
\end{pmatrix},
\end{equation}
couple the mode amplitudes. The superscripts on their entries indicate which types
of waves  they couple.  We have 
\begin{align}
\Gamma^{aa} (\bka,\bka',\zeta ) =&\frac{e^{ i k
 \nonumber 
   [\beta(\bka')-\beta(\bka)] \zeta}}{\beta^{\frac{1}{2}}(\bka)\beta^{\frac{1}{2}}(\bka')} , \quad \Gamma^{ab} (\bka,\bka',\zeta ) =\frac{e^{-i
  k [\beta(\bka')+\beta(\bka)  ] \zeta}}{\beta^{\frac{1}{2}}(\bka)\beta^{\frac{1}{2}}(\bka')}  ,
   \\ 
   \Gamma^{bb}
(\bka,\bka',\zeta ) =&- \overline{\Gamma^{aa} (\bka,\bka',\zeta )}, \quad 
  \Gamma^{ba} (\bka,\bka',\zeta ) = - \overline{\Gamma^{ab} (\bka,\bka',\zeta )},  \label{eq:Fzetaaa}
\end{align}
where the bar denotes complex conjugate, and substituting in (\ref{eq:twopointbvp}) we obtain the   energy conservation identity
\begin{equation*}
\int_{|\bka|<1} \frac{{\rm d}(k \bka)}{(2 \pi)^d} \Big[ |a^\ep(\bka,z)|^2  - |b^\ep(\bka,z)|^2 \Big] = ~\mbox{constant in } z.
\end{equation*}

\subsection{The Markov limit}
\label{sect:DL}
Here we describe the $\ep\to 0$ limit of the solution of 
(\ref{eq:twopointbvp}) with boundary conditions
(\ref{eq:WD15})-(\ref{eq:WD17}).  We begin in section \ref{sect:propag}
by writing the solution in terms of propagator matrices,
and   show in section \ref{sect:FScA} that we can neglect the backward and
evanescent waves. The limit of the forward going  amplitudes is
 in section \ref{sect:DiffA}.

\subsubsection{Propagator matrices}
\label{sect:propag}
The $2 \times 2$ propagator matrices ${\bf P}^\ep(\bka, z;\bka_o)$
are solutions of 
\begin{align}
\nonumber
 \partial_z {\bf P}^\ep(\bka,z;\bka_o)  \nonumber
=
\frac{i k \alpha}{2 \gamma^d \ep^{\frac{1}{2}}} \int_{|\bka'|<1} \frac{{\rm d}(k\bka')}{(2\pi)^d}
\widehat{\nu} \Big( \frac{k(\bka-\bka')}{\gamma}, 
\frac{\gamma z}{\ep}\Big) \\\times 
\boldsymbol{\Gamma}\Big(\bka,\bka',\frac{z}{\ep} \Big) {\bf 
  P}^\ep(\bka',z;\bka_o)  ,
\label{eq:ivp}
\end{align}
 for $z>0$, with initial condition ${\bf P}^\ep(\bka,z=0;\bka_o) = \delta(\bka-\bka_o) {\bf I}$, where ${\bf
  I}$ is the $2 \times 2$ identity matrix.  They allow us to write the solution of (\ref{eq:twopointbvp}) as 
\begin{equation}
\begin{pmatrix}
{a}^\ep(\bka,z) \\
b^\eps(\bka,z)
\end{pmatrix}
= \int_{|\bka_o|<1} {\rm d}\bka_o\, {\bf P}^\ep(\bka,z;\bka_o)
\begin{pmatrix}
{a}_o(\bka_o) \\
{b}^\ep(\bka_o,0)
\end{pmatrix},
\label{eq:ivp1}
\end{equation}
for all $ z>0$. In particular, when $z = L$, the backward going amplitude $b^\ep(\bka,L)$ in the left hand side vanishes by 
(\ref{eq:WD15}).

\subsubsection{The forward scattering approximation}
\label{sect:FScA}
Equation (\ref{eq:ivp1}) shows that the interaction of the forward and backward going wave amplitudes 
$a^\ep$ and $b^\eps$ depends on the coupling of the entries of the propagator. The $\ep\to 0$ limit of the propagator
\[{\bf P}^{\ep}(\bka,z;\bka_o) =
\begin{pmatrix} {P}^{aa,\ep}(\bka,z;\bka_o) & {P}^{ab,\ep}(\bka,z;\bka_o)\\ {P}^{ba,\ep}(\bka,z;\bka_o) & 
{P}^{bb,\ep}(\bka,z;\bka_o)
\end{pmatrix} 
\]
can be obtained and identified as a Markov process  that satisfies a system
of stochastic differential equations. We refer to   \cite{PK1974,PW1994} and 
 appendix \ref{ap:DL} for details. Here we  state  the results.

The stochastic differential equations for the limit entries of ${P}^{ab,\ep}(\bka,z;\bka_o)$ and ${P}^{ba,\ep}(\bka,z;\bka_o)$ 
are coupled to the limit entries of
${P}^{aa,\ep}(\bka',z;\bka_o)$ and ${P}^{bb,\ep}(\bka',z;\bka_o)$ through the coefficients
\[
\widetilde \cR\left( \frac{k(\vka-\vka^{' -})}{\gamma}\right) =  
\widetilde \cR \left( \frac{k(\bka-\bka')}{\gamma}, \frac{k(\beta(\bka)+\beta(\bka'))}{\gamma}\right),
\]
where  $\widetilde \cR$ is the power spectral density (\ref{eq:DefFT}) and 
$\vka = (\bka,\beta(\bka))$ and $\vka^{-} = (\bka,-\beta(\bka))$ are the wave vectors 
of the forward and backward going waves. 
The second argument in these coefficients comes from 
the phase factors
$\pm k (\beta(\bka)+\beta(\bka'))\zeta$  in the matrices $\Gamma^{ab}$ and $\Gamma^{ba}$.  
The coupling between  ${P}^{aa,\ep}(\bka,z;\bka_o)$ and $P^{aa,\eps}(\bka',z;\bka_o)$ is through the coefficients 
\[
\widetilde \cR\left( \frac{k(\vka-\vka')}{\gamma}\right) =  
\widetilde \cR \left( \frac{k(\bka-\bka')}{\gamma}, \frac{k(\beta(\bka)-\beta(\bka'))}{\gamma}\right),
\]
because the phase factors in matrices $\Gamma^{aa}$ are $ k (
\beta(\bka)-\beta(\bka') )\zeta $. The matrices $\Gamma^{bb}$ have
the same factors so the same coefficients 
couple the entries ${P}^{bb,\ep}$.

We conclude that  the coupling of the entries of the propagator  and therefore the interaction of the waves 
depends on the decay of the power 
spectral density $\widetilde \cR$.   
We now explain that when the mode amplitudes are supported initially at 
$|\bka| \le \eta/k <1$,   and $\gamma$ is as in \eqref{eq:F17},  we can  
neglect the backward going waves over distances of propagation of order $L$.

The power spectral density $\widetilde \cR( \vec\bq )$ is negligible
when $|\vec\bq|>1$, 
so  $\widetilde \cR( k \vec{\bka}/\gamma )$ is negligible
when $|\vec{\bka}| >  \gamma/k $. From (\ref{eq:sep}) and (\ref{eq:F17}), it is possible to choose some $\kappa_{_M} \in ({\eta} /k ,1)$ such that 
$\gamma$ satisfies
\begin{equation}
\label{eq:ASG2}
\frac{k \beta(\kappa_M)}{\gamma}  > 1.
\end{equation}
Then, for all $\bka'$ satisfying $|\bka'| < \kappa_{_M}$,
the coupling coefficients between $P^{aa,\ep}$ and $P^{ab,\ep}$ vanish
because 
\[
\frac{k |\vka-\vka^{'  -}|}{\gamma} \ge \frac{k(\beta(\bka) + \beta(\bka'))}{\gamma} \ge \frac{k \beta(\kappa_M)}{\gamma} > 1,
\]
and  $\widetilde \cR\big(k(\vka-\vka^{'-})/\gamma \big)$ is negligible.
  This implies the asymptotic decoupling 
of $a^\eps$ and $b^\eps$, and due to  the homogeneous boundary condition $b^\ep(\bka,L) = 0$,  
we conclude that we can neglect the backward going waves in the limit $\eps \to 0$.

The forward going amplitudes  interact  with each other,
because the coupling coefficients of the entries $P^{aa,\eps}$ of the propagator 
are large for at least a subset of transverse wave vectors satisfying $|\bka|,
|\bka'| \le \kappa_{_M}$ and 
\[
{|\bka-\bka'|},\, 
 {|\beta(\bka)-\beta(\bka')|}<  \frac{\gamma}{k}.
\]
Due to this coupling there is diffusion of
energy from the waves emitted by the source with $|\bka| < {\eta} / k $, to waves at larger values of $|\bka|$. This is why we take $\kappa_{_M}
> {\eta} / k$ in (\ref{eq:ASG2}). By assuming that
$a^\ep(\bka,z)$   are supported at $|\bka|
\le \kappa_{_M} < 1$ we essentially restrict $z$ by $Z_M$, so that the
energy does not diffuse to waves with $|\bka| > \kappa_{_M}$ for $z\le
Z_{_M}$. Physically,  the  wave vectors
$(\bka,\beta(\bka))$ of the forward going waves remain within a cone
with opening angle smaller than $180$ degrees.

We will see that the evolution of the $\bka$-distribution of the wave energy is described by a radiative transfer equation, 
which means that the wave energy undergoes a random walk (or diffusion). We can estimate 
from Eq. (\ref{eq:T7}) that the diffusion coefficient is of the order $\alpha^2 \gamma$, 
so the 
$\bka$-distribution of the wave energy 
reaches $\kappa_{_M}$ after a propagation distance of the order of $Z_{_M}$,  
such that $ \alpha^2 \gamma Z_{_M} = \kappa_{_M}^2$.
In dimensional units, this means $\alpha^2 Z_{_M} / \ell = \kappa_{_M}^2$.
Since $\alpha^2 L/ \ell =(\alpha^2 L/\lambda )(\lambda / \ell ) <1 $ by (\ref{eq:sep}),
it is possible to choose $Z_{_M} \sim L$ and a suitable $\kappa_{_M}<1$.

The evanescent waves can only couple with the propagating waves with
wave vectors of magnitude close to $1$. Thus, as long as the energy of
the wave is supported at $|\bka|< \kappa_{_M}$,
assumption (\ref{eq:ASG2}) implies that the evanescent waves do not
get excited.

\subsubsection{Markov limit of the forward going mode amplitudes}
\label{sect:DiffA}
We just explained that in the limit $\eps \to 0$ we can can neglect all the 
backward going waves and the evanescent ones. It remains to describe the 
limit of the forward going wave amplitudes $a^\ep(\bka,z)$ which satisfy the initial value
problem
\begin{align}
\nonumber
\partial_z a^\ep(\bka,z) =
  \frac{i k \alpha}{2 \gamma^d \ep^{\frac{1}{2}}} \int_{|\bka'| < 1} \frac{{\rm d} (k
    \bka')}{(2 \pi)^d}\, \widehat{\nu} \Big( \frac{k(\bka-\bka')}{\gamma}
  ,  \frac{\gamma z}{\ep}\Big)\\
  \times  \Gamma^{aa}\Big(\bka,\bka',\frac{z}{\ep} \Big)
a^\ep(\bka',z) ,
\label{eq:ivteps}
\end{align}
for $z > 0$, and the initial condition $a^\ep(\bka,0) =
  a_o (\bka) $.
These equations conserve energy, meaning that for all $\ep> 0$ and  all $z \geq0$,
\begin{equation}
\int_{|\bka| < 1} \frac{{\rm d} (k \bka)}{(2\pi)^d}|a^\ep(\bka,z)|^2 
= \int_{|\bka| < 1} \frac{{\rm d} (k \bka)}{(2\pi)^d}
 |a_o(\bka)|^2 \, .\label{eq:CENA}
\end{equation}

The details of the $\ep \to 0$ limit of $a^\ep(\bka,z)$ are in appendix \ref{ap:DL}.  
In particular, we explain there that the process 
\begin{align}
\bX^\ep(z) = \begin{pmatrix} {\rm Re} \big( {a}^\ep(\bka,z) )\\ {\rm
    Im}\big( {a}^\ep(\bka,z) ) 
\end{pmatrix}_{\hspace{-0.05in}  \bka \in  \mathscr{O} } ~~ \mbox{for} ~  
\mathscr{O} = \{\bka \in \RR^d, \, |\bka|<1\},
\label{def:Xeps}
\end{align}
converges weakly in ${\mathcal C}([0,L],{\mathcal D}')$ to a Markov process
$\bX(z)$, where ${\mathcal D}'$ is the space of distributions, dual to the
space $\mathcal{D}(\mathscr{O},\mathbb{R}^2)$ of infinitely
differentiable vector valued functions in $\mathbb{R}^2$, with compact
support. The generator of $\bX(z)$ is given in appendix \ref{ap:DL}, and we
denote henceforth the limit amplitudes by $(a(\bka,z))_{\bka \in \mathscr{O} }=X_1(z)+iX_2(z)$. Their first and second 
moments are described in the
next two sections.

\subsection{The coherent field}
\label{sect:COH}
The coherent wave field is 
\begin{align*}
  \EE\big[ u^\ep\big(\frac{\bx}{\eps},z\big)\big] \approx \int_{|\bka| <
    1} \frac{{\rm d}(k \bka)}{(2 \pi)^d} \frac{\EE[a(\bka,z)]}{
  \beta^{\frac{1}{2}}(\bka)}    e^{i \frac{k}{\ep}
    \vec{\bka} \cdot \vx} ,
\end{align*}
where  we replaced $\EE[a^\ep(\bka,z)]$ by its $\eps \to 0$ limit
$\EE[a(\bka,z)]$.
As explained in appendix 
\ref{ap:DL}, the mean field ${\mathcal A}(\bka,z) = \EE[a(\bka,z)]$ satisfies the
initial value problem
\begin{equation}
  \partial_z {\mathcal A}(\bka,z) = {\rm Q}(\bka)
{\mathcal A}(\bka,z), \quad z > 0,
  \label{eq:Amean}
\end{equation}
with initial condition ${\mathcal A}(\bka,0) =
a_o(\bka)$, and ${\rm Q}(\bka)$ given by
\begin{align}
{\rm Q}&(\bka) = -\frac{k^2 \alpha^2}{4 \gamma^{d+1}} \int_{|\bka'| < 1} \frac{{\rm d}(k \bka')}{(2 \pi)^d} 
\frac{1}{\beta(\bka)\beta(\bka')} \nonumber \\&
\times \int_0^\infty {\rm d} \zeta \int_{\mathbb{R}^d} {\rm d} \br \, \cR(\br,\zeta) e^{- i \frac{k}{\gamma} 
\big(\bka-\bka',\beta(\bka)-\beta(\bka')\big) \cdot (\br,\zeta)}.
\label{eq:defQ}
\end{align}
This is the same as \eqref{eq:C5} in our scaling. 

The solution of (\ref{eq:Amean}) is 
\begin{equation}
{\mathcal A}(\bka,z)= \exp\big[ {\rm Q}(\bka)z \big]
a_o(\bka),
\label{eq:C6}
\end{equation}
so as {stated} in section \ref{sect:MR}, the random medium effects do not average
out. The mean amplitudes are not the same as the amplitudes in the
homogeneous medium at $z > 0$, and they decay with $z$ on the $\bka$ dependent scales ${\mathcal S}(\bka) = - 1/{{\rm Re}[{\rm Q}(\bka)]} ,$ 
the scattering mean free paths. The real part 
of ${\rm Q}(\bka)$, which is non-positive,  is
an effective diffusion term in \eqref{eq:Amean}, which removes energy
from the mean field and gives it to the incoherent fluctuations. This is due to 
the randomization or loss of coherence of the waves. 
 The
imaginary part of $ Q(\bka)$ 
is an effective dispersion term, which does not remove energy from the
mean field and ensures causality\footnote{ If we write the coherent
  wave fields in the time domain, using the inverse Fourier transform
  with respect to the frequency $\om$, we obtain a causal result.}.

%
\subsection{The one-way radiative transfer equations}
\label{sect:TEQ}
The mean intensity in the direction of $\bka$ is 
\begin{equation}
 {{\mathcal I}}(\bka,z) = \lim_{\ep\to 0} \EE \left[ | a^\ep(\bka,z) |^2 \right],
  \label{eq:T1}
\end{equation}
and it evolves in $z> 0$ as modeled by equation 
\begin{align}
 \partial_z   {{\mathcal I}}(\bka,z) =& \int_{|\bka'| < 1} \frac{{\rm d} (k \bka')}{(2 \pi)^d}\,
{\mathcal Q}(\bka,\bka')  \big[  {{\mathcal I}}(\bka',z)  -  {{\mathcal I}}(\bka,z) \big],
  \label{eq:T2b}
\end{align}
with initial condition $  {{\mathcal I}}(\bka,0) =|a_o(\bka)|^2$ (see Appendix \ref{ap:DL}).
The differential scattering cross section
\begin{align*}
\hspace{-0.09in} {\mathcal Q}(\bka,\bka') = \frac{k^2 \alpha^2}{4 \gamma^{d+1}\beta(\bka)\beta(\bka') } 
\widetilde \cR \Big( \frac{k}{\gamma} \big(\bka-\bka',\beta(\bka)-\beta(\bka')\big) \Big)
\end{align*}
is the same as (\ref{eq:DCS}) in our scaling, and from (\ref{eq:defQ}) we see that $-2{\rm Re}[{\rm Q}(\bka)]$ equals the 
total scattering cross section
\begin{equation}
-2{\rm Re}[{\rm Q}(\bka)]=  \int_{|\bka'| < 1} \frac{{\rm d} (k \bka')}{(2 \pi)^d}\,
{\mathcal Q}(\bka,\bka').
\label{eq:TSC}
\end{equation}
We also note that  the intensities satisfy the 
conservation identity
\begin{equation*}
  \int_{|\bka|<1} \frac{{\rm d} (k \bka)}{(2\pi)^d}  
    {{\mathcal I}}(\bka,z)  = \int_{|\bka|<1} \frac{{\rm d} (k
    \bka)}{(2\pi)^d}  |a_o(\bka)|^2, \quad \mbox{for all } z > 0,  
\end{equation*}
which is consistent with (\ref{eq:CENA}).

Using the generator of the Markov limit process ${\bf X}(z)$ given in appendix \ref{ap:DL}, 
we can also calculate the $\ep \to 0$ limit of the second moments $\EE \left[  
      a^\ep(\bka,z)\overline{a^\ep(\bka',z) } \right] $ of the mode amplitudes. We obtain that 
when $\bka \ne \bka'$, 
\begin{align*}
 \lim_{\ep\to 0} \EE \left[  
      a^\ep(\bka,z)\overline{a^\ep(\bka',z) } \right] 
      = \lim_{\ep\to 0}  
      \EE[a^\ep(\bka,z)] 
\overline{   \EE[a^\ep(\bka',z)] }, 
\end{align*}
meaning that the waves traveling in different directions are asymptotically decorrelated \footnote{It can also be shown that the waves decorrelate over
frequency offsets larger than $\ep$. Thus, one can study the energy
density resolved over both time and space i.e., the space-time Wigner
transform. 
}. This is because  these waves  see  different regions of the random medium. It is only when the waves propagate in 
similar directions i.e., $|\bka'-\bka| = O(\ep)$, that the mode amplitudes are correlated, so we define the energy density (Wigner transform) as 
\begin{align}
{\mathcal W}(\bka,\bx,z) = \lim_{\ep\to 0} \int \frac{{\rm d} (k {\itbf q})}{(2 \pi)^d}
\, \exp\Big[ i k {\itbf q} \cdot \big(\nabla \beta(\bka) z + \bx\big)
  \Big] 
  ~~\nonumber \\ \times 
  \EE \left[   a^\ep\big(\bka
    + \frac{\ep{\itbf q}}{2},z\big) \overline{    a^\ep\big(\bka-\frac{\ep{\itbf q}}{2},z\big)}\right].
\label{eq:T6}
\end{align}
It satisfies the transport equation 
\begin{align}
\nonumber
& \partial_z {\mathcal W}(\bka,\bx,z) - \nabla \beta(\bka) \cdot \nabla_\bx
  {\mathcal W}(\bka, \bx,z)\\
  & = \int_{|\bka'| < 1} \frac{{\rm d} (k \bka')}{(2 \pi)^d}\,
{\mathcal Q}(\bka,\bka')  \big[  {\mathcal W}(\bka',\bx,z)  - {\mathcal W}(\bka,\bx,z)  \big] ,
  \label{eq:T7}
\end{align}
for $z>0$, as stated in section \ref{sect:MR}. 
When the initial condition $a_o(\bka)$ is smooth in $\bka$, we have from \eqref{eq:T6} that 
\begin{equation*}
{\mathcal W}(\bka,\bx,0) = \delta(\bx)
|a_o(\bka)|^2,
\end{equation*}
and therefore at $z >0$
\begin{equation*}
{\mathcal W}(\bka,\bx,z) = \delta\big(\bx + \nabla \beta(\bka) z\big)
{{\mathcal I}}(\bka ,z ).
\end{equation*}
This shows that the energy is transported on the characteristic 
\[
\bx = -\nabla \beta(\bka) z = \frac{\bka}{\beta(\bka)} z.
\]

\section{The high-frequency limit}
\label{sect:HF}
In the high-frequency limit $\gamma \to 0$ the transport equations simplify. We
quantify the scattering mean free paths in this limit, and show how to derive the diffusion approximation
and  paraxial model from the transport equations \eqref{eq:T7}.

\subsection{Quantification of scattering mean free paths}
\label{sect:HF1}
If we expand in powers of $\gamma$ the right hand side of \eqref{eq:TSC},
we obtain the following expression of the scattering mean free paths
\begin{align*} 
\cS(\bka) &= - \frac{1}{{\rm Re} \big[ {\rm Q}(\bka)\big]} =  \frac{8  \gamma \beta^2(\bka)}{k^2 \alpha^2
   \int_{-\infty}^\infty  {\rm d} \zeta\, \cR\Big(\frac{\bka \zeta}{\beta(\bka)},\zeta\Big)}
+ O(\gamma^2).
\end{align*}
They are of order $\gamma$ and decrease as the negative power of $2$
with the frequency $\omega = k c_o$, meaning that higher frequency waves lose coherence faster.
We also expect that $\cS(\bka)$ decrease monotonically with $|\bka|$, because a plane wave mode with wavevector $k(\bka,\beta(\bka))$  travels the distance $z/\beta(\bka)$ in the random medium when it 
propagates up to  $z$. The closer $|\bka|$ is to one, the longer the distance and thus, the faster 
the loss of coherence quantified by the scale $\cS(\bka)$.  The monotone dependence of $\cS(\bka)$ on 
$|\bka|$  can be seen explicitly in statistically isotropic media, where $\cR(\vec\bx)=\cR_{\rm iso}(|\vec\bx|)$, and 
\[
\cR \Big(\frac{\bka \zeta}{\beta(\bka)},\zeta \Big) = \cR_{\rm iso} \left( \sqrt{\frac{|\bka|^2 \zeta^2}{\beta^2(\bka)} + \zeta^2} \right) = 
\cR_{\rm iso}\Big(\frac{|\zeta|}{\beta(\bka)}\Big).
\]
Then
  \begin{align*}
\cS(\bka) = \frac{4 \gamma 
  \beta(\bka)}{ k^2 \alpha^2 \int_0^\infty   {\rm d} \zeta \,  \cR_{\rm iso} (\zeta)} + O(\gamma^2),
\end{align*}
and the decay with $|\bka|$ is captured by $\beta(\bka) = \sqrt{1-|\bka|^2}$.

\subsection{The diffusion approximation}
\label{sect:HF2}
The mean mode intensities ${\mathcal I}(\bka,z)$ defined in (\ref{eq:T1}) satisfy  (\ref{eq:T2b}), 
with initial condition at $z=0$ derived from (\ref{eq:WD8bp}):
\[
{\mathcal I}(\bka,0) = \frac{1}{4 k^2 \beta(\bka) \eta^{2d}} \Big| \what F\Big(\frac{k \bka}{\eta} \Big) \Big|^2.
\]
This is independent of $\gamma$ and for fixed $\eta$.
 
The diffusion model is obtained by expanding Eq.~(\ref{eq:T2b}) in powers of $\gamma$. We obtain that 
\begin{align}
\partial_z {\mathcal I}(\bka,z) \approx \gamma
\bigg[\sum_{j,l=1}^d  A_{jl}(\bka)  \partial^2_{\kappa_j \kappa_l} 
+\gamma
\sum_{j=1}^d B_j(\bka)  \partial_{\kappa_j }  \bigg]
{\mathcal I}(\bka,z),
\label{eq:HC16b}
\end{align}
where the approximation means that we neglect higher powers  in $\gamma$, and the 
diffusion and drift coefficients are independent of $k$ and $\gamma$:
\begin{align*}
A_{jl}(\bka) =& -  \frac{\alpha^2}{8\beta(\bka)^2 }\int_{-\infty}^\infty  {\rm d} \zeta \, \partial^2_{r_jr_l}   \cR\Big(\frac{\bka \zeta}{\beta(\bka)},\zeta\Big) , 
\quad \quad j,l=1,\ldots,d,
\end{align*}
and 
\begin{align*}
B_j(\bka) =& \sum_{l,m=1}^d
\frac{\alpha^2 \partial^2_{\kappa_l\kappa_{_M}} \beta(\bka) }
{8\beta(\bka)^2 }
\int_{-\infty}^\infty  {\rm d} \zeta \, \zeta \partial^3_{r_jr_l r_m}   \cR\Big(\frac{\bka \zeta}{\beta(\bka)},\zeta\Big)  \\
&
- \sum_{l=1}^d
\frac{\alpha^2 \kappa_l}{4\beta(\bka)^4 } 
\int_{-\infty}^\infty  {\rm d} \zeta \,   \partial^2_{r_jr_l}   \cR\Big(\frac{\bka \zeta}{\beta(\bka)},\zeta\Big) 
 ,\quad \quad j=1,\ldots,d .
\end{align*}
Note that the diffusion is the dominant term in (\ref{eq:HC16b}).

\subsection{The paraxial approximation}
The paraxial (beam-like) propagation model is for a large diameter $X$ of the support of the source 
with respect to the wavelength, so that  $\eta \to 0$. The result depends on the order in which we take the limits $\eta \to 0$ and $\gamma \to 0$, as we now explain.

In regimes with $\la \ll \ell  = X$, where $\eta = \gamma$,  the rescaled intensity
\begin{align*}
{\mathcal I}_{\rm res}(\bka,z) = \gamma^{2d} {\mathcal I}\big(  \gamma \bka, \gamma z\big)  
\end{align*}
satisfies  in the limit  $\gamma \to 0$ the equation
\begin{align}
\label{eq:rteparlim}
\hspace*{-0.1in}
\partial_z  {\mathcal I}_{\rm res}= \frac{k^2 \alpha^2}{4} \int_{\RR^d} \frac{{\rm d} (k\bka')}{(2\pi)^d} \widetilde{\cR}\big(k(\bka-\bka'),0\big)
\big[ {\mathcal I}_{\rm res}(\bka') - {\mathcal I}_{\rm res} (\bka)  \big] ,
\end{align}
with initial condition 
$
{\mathcal I}_{\rm res}(\bka,0) = \big| \what F({k \bka}) \big|^2 / [4 k^2 \beta(\bka)]
$.
This is the transport equation for the random paraxial wave  equation, as explained in subsection \ref{app:C}.

In regimes with $\la \ll \ell \ll X$,  analyzed with the sequence of limits $\gamma \to 0$, followed by $\eta \to 0$, the rescaled intensity
\begin{align*}
{\mathcal I}_{\rm res}(\bka,z) = \eta^{2d} {\mathcal I}\Big( {\eta} \bka, \frac{{\eta}^2}{\gamma} z\Big)  
\end{align*}
satisfies  the diffusion equation 
\begin{align}
\partial_z  {\mathcal I}_{\rm res} = \sum_{j,l=1}^d
 {D}_{{\rm res},jl} \partial^2_{\kappa_j\kappa_l} {\mathcal I}_{\rm res},
 \label{eq:difpar}
\end{align}
with initial condition
$
{\mathcal I}_{\rm res}(\bka,0) = \big| \what F({k \bka}) \big|^2 / [4 k^2 \beta(\bka)]
$
and diffusion tensor ${D}_{{\rm res}, jl}  $ given by 
\begin{align*}
{D}_{{\rm res},jl}
=
- \frac{\alpha^2}{8}
\int_{-\infty}^\infty {\rm d}\zeta \,  \partial_{r_jr_l}^2  \cR ({\bf 0},\zeta )  = \lim_{|\bka| \to 0} A_{jl}(\bka),
\end{align*}
for $j,l=1,\ldots,d$.
This result was derived in \cite{Fannjiang,papa07,garnier2009paraxialderiv,garnier2009paraxial} starting from the 
paraxial wave equation. We recovered it here because in the regime  with $\la \ll \ell \ll X$ we have 
a narrow cone beam propagating through a random medium.

Note that equation (\ref{eq:difpar}) can also be derived formally  from
the radiative transfer equation (\ref{eq:RTT}). 
First, one considers that scattering is sharply peaked in the forward scattering direction,
so that it is possible to take the Fokker-Planck approximation,
that is to say, the right-hand side of (\ref{eq:RTT}) can be approximated by a diffusion operator in $\vk$ \cite{pomraning,larsen}.
Second, one considers that the source emission is sharply peaked and that the propagation distance 
is short enough so that the wave remains in the form of a narrow cone beam.

\section{Summary}
\label{sect:sum}
The one-way radiative transfer equation describes the evolution of the intensity of the waves 
resolved over directions, the Wigner transform, in forward-peaked scattering regimes. We  derived it
using multiscale analysis and probabilistic limits, starting from the wave equation in random 
media. The scattering regime with small random fluctuations of the wave speed and long distances
of propagation over which cumulative scattering becomes significant leads to waves 
propagating forward in a wide angular cone. It bridges between two known regimes: The first is the radiative
transfer regime where waves propagate in all directions and the Wigner transform satisfies 
a boundary value problem. The second is the paraxial regime, where waves propagate in a  narrow angle 
cone. We established this bridge by connecting  the one-way radiative transfer equation to 
the equations for the Wigner transform in these two regimes.

\section*{Acknowledgements}
We thank John Schotland for useful discussions and encouragement to pursue this project.
Liliana Borcea's work was partially supported by the NSF grant DMS1510429. Support from
the Simons Foundation is also gratefully acknowledged.

\appendix
\section{The Markov limit}
\label{ap:DL}
Let $\mathscr{O}$ be an open set in $\RR^d$ and ${\mathcal D}(\mathscr{O},\RR^2)$ the space of
infinitely differentiable functions with compact support.  We consider
the process $\bX^\ep$ in the space   ${\mathcal C}([0,L], {\mathcal D}'(\mathscr{O},\RR^2))$ of continuous functions of $z$. It is  the solution of
\begin{equation}
\frac{{\rm d} \bX^\ep}{{\rm d}z} = \frac{1}{\sqrt{\ep}} {\mathcal
  F}\Big(\frac{z}{\ep},\frac{z}{\ep}\Big) \bX^\ep ,
  \label{eq:generalformapp}
\end{equation}
where ${\mathcal F}(\zeta,\zeta')$ is a random
linear operator from ${\mathcal D}'$ to ${\mathcal D}'$.  Here ${\mathcal D}'$ denotes the space of 
distributions, dual to ${\mathcal D}(\mathscr{O},\RR^2)$.
We assume that the
mapping $\zeta \to {\mathcal F}(\zeta,\zeta')$ is stationary and possesses strong ergodic
properties, and that ${\mathcal F}(\zeta,\zeta')$ has mean zero. Moreover,
the mapping $\zeta' \to {\mathcal F}(\zeta,\zeta')$ is periodic.

We are interested in particular in equation (\ref{eq:ivteps}), that can be put into the form (\ref{eq:generalformapp}) if we define 
the process $\bX^\ep$  as (\ref{def:Xeps})
and the operator ${\mathcal F}(\zeta,\zeta')$ as
\begin{align}
\left< {\mathcal F}(\zeta,\zeta') \bX ,\bphi\right> = \sum_{j=1}^2 \int_\mathscr{O} \frac{{\rm d}(k\bka)}{(2\pi)^d}
 [{\mathcal F}(\zeta,\zeta') \bX ]_j (\bka) \phi_j(\bka) \nonumber \\ 
= \int_\mathscr{O} \frac{{\rm d}(k\bka)}{(2\pi)^d} \bphi(\bka) \cdot \int_\mathscr{O} \frac{{\rm d}(k\bka')}{(2\pi)^d}
\mathbfcal{F} (\bka,\bka',\zeta,\zeta') {\itbf X}(\bka'),
\label{eq:A3}
\end{align}
for $\boldsymbol{\phi} \in {\mathcal D}(\mathscr{O},\RR^2)$ with components
$\phi_j$ and  ${\itbf X} \in {\mathcal D}'(\mathscr{O},\RR^2)$ with components
$X_j$. The kernel matrix $\mathbfcal{F}(\bka,\bka',\zeta,\zeta')$ is given by
\begin{align}
\mathbfcal{F}= \begin{pmatrix} {\mathscr F}^r &
  -{\mathscr F}^i 
  \\ {\mathscr F}^i & {\mathscr F}^r   \end{pmatrix},
\label{eq:A4}
\end{align}
in terms of
\begin{align}
{\mathscr F}^{\rm r}(\bka,\bka',\zeta,\zeta') &= {\rm Re}
\Big[\frac{i k \alpha}{2 \gamma^d} \widehat{\nu} \Big(
\frac{k(\bka-\bka')}{\gamma},\gamma \zeta \Big) \Gamma^{aa}
  \big(\bka,\bka',\zeta' \big) \Big]
,\label{eq:A5}\\ {\mathscr F}^{\rm i}(\bka,\bka',\zeta,\zeta') &=
    {\rm Im} \Big[\frac{i k \alpha}{2\gamma^d} \widehat{\nu}
      \Big( \frac{k(\bka-\bka')}{\gamma},\gamma \zeta\Big) \Gamma^{aa}
      \big(\bka,\bka',\zeta' \big) \Big] , \label{eq:A6}
\end{align}
where we recall from (\ref{eq:Fzetaaa}) the
expression of $\Gamma^{aa}(\bka,\bka',\zeta')$.  The adjoint operator
${\mathcal F}^*(\zeta,\zeta')$ is defined by
$$ \left< {\mathcal F}(\zeta,\zeta') \bX ,\bphi\right> = \left< \bX ,{\mathcal
  F}^*(\zeta,\zeta') \bphi\right>
$$ for $\bphi \in {\mathcal D}(\mathscr{O},\RR^2)$ and $\bX \in {\mathcal
  D}'(\mathscr{O},\RR^2)$, and has matrix kernel
$\mathbfcal{F}^*(\bka,\bka',\zeta,\zeta') = \mathbfcal{F}
(\bka',\bka,\zeta,\zeta')^T$, where the superscript $T$ stands for transpose.

To obtain the Markov limit we use the results in \cite{PW1994} (the interested reader may first read 
\cite[Chap. 6]{fouque07} for a self-contained introduction to such limit theorems).
They give that $\bX^\ep(z)$ converges weakly in ${\mathcal C}([0,L],{\mathcal
  D}')$ to $\bX(z)$,  the solution of a martingale problem
with generator ${\mathcal L}$ defined by
\begin{align}
\nonumber
& {\mathcal L} f(\left<\bX , \bphi \right>) =\\
\nonumber
& \int_0^\infty  {\rm d} \zeta \lim_{Z\to \infty} \frac{1}{Z}
\int_0^Z {\rm d} h \, \EE \big[ \left< \bX , {\mathcal F}^*(0,h)
  \bphi\right>   \left< \bX , {\mathcal F}^*(\zeta,\zeta+h) \bphi\right>
  \big] \\
  \nonumber
  &\quad \times f'' ( \left< \bX , \bphi\right> ) \\ 
  \nonumber
&+
\int_0^\infty {\rm d} \zeta \lim_{Z\to \infty} \frac{1}{Z}
\int_0^Z {\rm d} h \, \EE \big[ \left< \bX , {\mathcal F}^*(0,h) 
  {\mathcal F}^*(\zeta,\zeta+h) \bphi\right> \big] \\
& \quad \times  f' ( \left< \bX ,
\bphi\right> )  ,
  \label{eq:GEN0}
\end{align}
for any $\bX \in {\mathcal D}'(\mathscr{O},\RR^2)$, $\bphi\in {\mathcal
  D}(\mathscr{O},\RR^2)$, and smooth $f:\RR\to \RR$.  This
means that, for any $\bphi \in {\mathcal D}(\mathscr{O},\RR^2)$ and smooth
function $f :\RR\to \RR$, the real-valued process
$$
 f \big(\left< \bX(z) , \bphi \right> \big)- \int_0^z {\rm d} z' \, {\mathcal L} f \big(\left<
\bX(z') , \bphi \right> \big)
$$
 is a martingale.  More generally, if $n \in \NN$,
$\bphi^{(1)},\ldots,\bphi^{(n)}\in {\mathcal D}(\mathscr{O},\RR^2)$, and $f:\RR^n \to \RR $
is a smooth function, then
\begin{align}
&f \big(\left< \bX(z) , \bphi^{(1)}\right>,\ldots,\left< \bX(z) ,
\bphi^{(n)}\right>\big) \nonumber\\
&- \int_0^z {\rm d} z' \, {\mathcal L}^{(n)} f \big( \left< \bX(z') ,
\bphi^{(1)} \right>,\ldots,\left< \bX(z') , \bphi^{(n)} \right> \big) 
\label{eq:Mart}
\end{align}
is a martingale,
where 
\begin{align}
&  {\mathcal L}^{(n)} f \big(\left< \bX , \bphi^{(1)} \right>,\ldots,\left< \bX ,
  \bphi^{(n)}  \right>\big)\nonumber =\\
\nonumber
&   \sum_{j,l=1}^n
  \int_0^\infty {\rm d} \zeta \lim_{Z\to \infty}
  \frac{1}{Z}\int_0^Z {\rm d} h \, \EE \big[ \left< \bX , {\mathcal F}^*(0,h) \bphi^{(j)} \right>   \\
    &\quad \times  \left< \bX , {\mathcal
      F}^*(\zeta,\zeta+h) \bphi^{(l)} \right> \big]\, \partial^2_{jl} f(\left< \bX ,
  \bphi^{(1)}  \right>,\ldots,\left< \bX , \bphi^{(n)} \right>\big)\nonumber \\ 
 \nonumber &
+ \sum_{j=1}^n
  \int_0^\infty  {\rm d} \zeta \lim_{Z\to \infty}
  \frac{1}{Z}\int_0^Z {\rm d} h \, \EE \big[ \left< \bX ,
    {\mathcal F}^*(0,h) {\mathcal F}^*(\zeta,\zeta+h) \bphi^{(j)} \right> \big] \\
    &\quad \times 
  \partial_{j} f(\left< \bX ,
  \bphi^{(1)}  \right>,\ldots,\left< \bX , \bphi^{(n)}  \right>\big) . \label{eq:GEN}
\end{align}

To calculate the first moment of the limit process $\bX(z)$, let $n =
1$ and $f(y)=y$ in (\ref{eq:Mart})-(\ref{eq:GEN}). We find that
\begin{eqnarray*}
\frac{ {\rm d} \, \EE \big[ \left< {\bX}(z) , \bphi \right>\big]}{{\rm d} z} = \EE
\big[ \left< \bX(z) , {\mathcal H}^*\bphi \right>\big]   ,
\end{eqnarray*}
where 
\begin{eqnarray*}
{\mathcal H}^* =\int_0^\infty {\rm d} \zeta
\lim_{Z\to \infty}\frac{1}{Z}\int_0^Z {\rm d} h \, \EE
\big[ {\mathcal F}^*(0,h) {\mathcal F}^*(\zeta,\zeta+h) \big]
 .
\end{eqnarray*}
This shows that
$$
\boldsymbol{\mathcal X}(z) = \EE \big[ \bX(z) \big] 
$$ 
satisfies a closed system of ordinary differential equations
\begin{eqnarray*}
\frac{ {\rm d} \left< \boldsymbol{\mathcal X}(z) , \bphi \right> }{{\rm d} z} = \left<
\boldsymbol{\mathcal X}(z) , {\mathcal H}^*\bphi \right>   ,
\end{eqnarray*}
or, equivalently in ${\mathcal D}'$,
\begin{equation}
\frac{ {\rm d} \boldsymbol{\mathcal X}(z) }{{\rm d} z} = {\mathcal H}
\boldsymbol{\mathcal X}(z)   ,
\label{eq:MeanField}
 \end{equation}
where ${\mathcal H}$ is the adjoint
of ${\mathcal H}^*$. The kernel matrix of ${\mathcal H}$ is $\mathbfcal{H}(\bka',\bka) = \mathbfcal{H}^*(\bka,\bka')^T$.
Recalling from (\ref{eq:A3})-(\ref{eq:A6}) the expression of the
kernel $\mathbfcal{F}(\bka',\bka,\zeta,\zeta')^T$ of ${\mathcal F}^*(\zeta,\zeta')$, we obtain
that the matrix kernel  $\mathbfcal{H}^*(\bka',\bka)$ of ${\mathcal H}^*$ is
\begin{align*}
{\mathcal H}^*_{jl}(\bka,\bka') = \sum_{q=1}^2 \int_\mathscr{O} \frac{{\rm d}(k\bka'')}{(2\pi)^d}
\int_0^\infty  {\rm d} \zeta \lim_{Z\to \infty} \frac{1}{Z}
\int_0^Z   {\rm d} h  \\
\times 
\EE \big[{\mathcal F}_{lq}(\bka',\bka'',\zeta,\zeta+h) {\mathcal F}_{qj}(\bka'',\bka,0,h) \big] ,
\end{align*}
for $j,l = 1, 2$.  For instance,
\begin{align*}
{\mathcal H}^*_{11}(\bka,\bka') =& \int_\mathscr{O} \frac{{\rm d}(k\bka'')}{(2\pi)^d}\int_0^\infty
  {\rm d} \zeta \lim_{Z\to \infty} \frac{1}{Z}
\int_0^Z  {\rm d} h \\& \quad \times  \EE \big[ {\mathscr F}^{\rm
    r} (\bka',\bka'',\zeta,\zeta+h) {\mathscr F}^{\rm
    r} (\bka'',\bka,0,h) \big]\\ &-
 \int_\mathscr{O}\frac{{\rm d}(k\bka'')}{(2\pi)^d}\int_0^\infty
 {\rm d} \zeta \lim_{Z\to \infty} \frac{1}{Z}
\int_0^Z  {\rm d} h \\
& \quad \times\EE \big[
  {\mathscr F}^{\rm i} (\bka',\bka'',\zeta,\zeta+h) {\mathscr
    F}^{\rm i} (\bka'',\bka,0,h) \big]   ,
\end{align*}
and using (\ref{eq:A5})-(\ref{eq:A6}), we get 
\begin{align*}
{\mathcal H}^*_{11}(\bka,\bka') =& {\rm Re} \Big\{ \Big(
\frac{ik \alpha}{2 \gamma^d}\Big)^2\int_\mathscr{O} \frac{{\rm d}(k\bka'')}{(2\pi)^d}
\int_0^\infty   {\rm d} \zeta\lim_{Z\to \infty} \frac{1}{Z}
\int_0^Z  {\rm d} h \, \\ & \times \EE\Big[
  \widehat{\nu}\Big(\frac{k(\bka'-\bka'')}{\gamma},
  \gamma\zeta\Big)\widehat{\nu}\Big(\frac{k(\bka''-\bka)}{\gamma},0\Big)
  \Big] \\
 & \times \big[ \Gamma^{aa}
  (\bka',\bka'',\zeta+h) \Gamma^{aa}(\bka'',\bka,h) \big]
\Big\} .
\end{align*}
Moreover, using the identity
\begin{align*}
&\EE \Big[
  \widehat{\nu}\Big(\frac{k(\bka'-\bka'')}{\gamma},\gamma\zeta\Big)
\widehat{\nu}\Big(\frac{k(\bka''-\bka)}{\gamma},0\Big)\Big]
= \\
&\big(  2\pi \gamma \big)^d \delta\big( k (\bka-\bka') \big) \widehat{\cR}\Big( \frac{k(\bka-\bka'')}{\gamma},\gamma\zeta \Big) ,
 \end{align*}
 with
 \begin{align*}
 \what \cR ({\itbf q},\zeta) = \int_{\RR^d} \cR(\br,\zeta) e^{-i {\itbf q} \cdot \br} {\rm d} \br ,
 \end{align*}
derived from the definition of the autocovariance function with
straightforward algebraic manipulations, and obtaining from
(\ref{eq:Fzetaaa}) that 
\begin{eqnarray*}
\Gamma^{aa} (\bka,\bka'',\zeta+h) \Gamma^{aa}(\bka'',\bka,h) =
\frac{1}{\beta(\bka)\beta(\bka'')}e^{i k
  (\beta(\bka'')-\beta(\bka)) \zeta},
\end{eqnarray*}
we get 
\begin{align*}
{\mathcal H}^*_{11}(\bka,\bka') = -\frac{k^{2} \alpha^2}{ 4 \gamma^d}
         {\rm Re}\Big\{ \int_\mathscr{O} \frac{{\rm d}(k\bka'')}{(2\pi)^d} \int_0^\infty {\rm d} \zeta \,
           \widehat{\cR}\Big(
           \frac{k(\bka-\bka'')}{\gamma},\gamma\zeta \Big) \\
           \times e^{i k
             (\beta(\bka'')-\beta(\bka)) \zeta}  \frac{(2\pi)^d }{\beta(\bka)\beta(\bka'')} \,\delta\big( k(\bka-\bka') \big) \Big\}.
\end{align*}
The expressions of the other components of ${\mathcal H}^*_{jl}(\bka,\bka')$ are of the same type. Substituting into
(\ref{eq:MeanField}) we obtain the explicit expression of the
differential equations satisfied by the mean wave amplitudes. This is
equation (\ref{eq:Amean}), written in complex form.

The calculation of the second moments is similar, by letting $n = 1$
and $f(y) = y^2$ in (\ref{eq:GEN}), and carrying the lengthy
calculations.

\bibliographystyle{siam}
\bibliography{SCALAR}

\end{document}